\documentclass[12]{amsart}
\usepackage{amsmath,amssymb,amsthm,color,enumerate,comment,centernot,enumitem,url,cite}
\usepackage{graphicx,relsize,bm}
\usepackage{mathtools}
\usepackage{array}

\makeatletter
\newcommand{\tpmod}[1]{{\@displayfalse\pmod{#1}}}
\makeatother

\newtheorem{thm}{Theorem}[section]
\newtheorem{lemma}[thm]{Lemma}

\theoremstyle{remark}

\theoremstyle{definition}

\newtheorem{rem}[thm]{Remark}

\theoremstyle{THM}

\newcommand{\G}{{\mathcal G}}

\newcommand{\FF}{{\mathcal F}}

\newcommand{\mmod}[1]{\ \mathrm{mod}\enspace #1}
\newcommand{\Z}{{\mathbb Z}}
\newcommand{\Q}{{\mathbb Q}}

\newcommand{\F}{{\mathbb F}}
\newcommand{\abs}[1]{\left|{#1}\right|}

\newcommand{\ds}{\displaystyle}
\newcommand{\D}{{\mathcal D}}

\makeatletter
\@namedef{subjclassname@2020}{%
  \textup{2020} Mathematics Subject Classification}
\makeatother

\def\red#1 {\textcolor{red}{#1 }}
\def\blue#1 {\textcolor{blue}{#1 }}

\numberwithin{equation}{section}

\begin{document}

\title[The monogenicity of $\ds x^n+c\sum_{j=0}^{n-1}a^{n-j-1}x^j$]
{The monogenicity of $\ds \boldsymbol{x^n+c\sum_{j=0}^{n-1}a^{n-j-1}x^j}$}

\author{Lenny Jones}
\address{Professor Emeritus, Department of Mathematics, Shippensburg University, Shippensburg, Pennsylvania 17257, USA}
\email[Lenny~Jones]{doctorlennyjones@gmail.com\\ ORCID: 0000-0001-7661-4226}

\date{\today}

\begin{abstract}
We say that a monic polynomial $f(x)\in \Z[x]$ is {\em monogenic} if $f(x)$ is irreducible over ${\mathbb Q}$ and $\{1,\theta,\theta^2,\ldots,\theta^{\deg(f)-1}\}$ is a basis for the ring of integers of ${\mathbb Q}(\theta)$, where $f(\theta)=0$. Let $n,a,c\in \Z$ with $n\ge 3$, and suppose that 
\[\FF_{n,a,c}(x):=x^n+cx^{n-1}+acx^{n-2}+a^2cx^{n-3}+\cdots +a^{n-2}cx+a^{n-1}c\] is irreducible over ${\mathbb Q}$. 
Recently, Jakhar, Kalwaniya and Yadav have given necessary and sufficient conditions for the monogenicity of $\FF_{n,a,c}(x)$, under the assumption that all prime divisors of the discriminant of $\FF_{n,a,c}(x)$ not dividing $ac$ also do not divide $n+1$. 
In this note, we provide necessary and sufficient conditions for the monogenicity of $\FF_{n,a,c}(x)$ without this restriction. 
\end{abstract}


\maketitle
\section{Introduction}\label{Section:Intro}

Let $f(x)\in \Z[x]$ be monic and irreducible over $\Q$ with $f(\theta)=0$. Let $K=\Q(\theta)$, and let $\Z_K$ denote the ring of integers of $K$.  Then  
\begin{equation}\label{Eq:Dis-Dis}
\Delta(f)=\left[\Z_K:\Z[\theta]\right]^2\Delta(K), \quad \text{\cite{Cohen}}
\end{equation} where $\Delta(f)$ and $\Delta(K)$ denote the discriminants over $\Q$, respectively, of $f(x)$ and the number field $K$. We refer to $\left[\Z_K:\Z[\theta]\right]$ as the {\em index} of $f(x)$, and we denote it index($f$). Furthermore, we say $f(x)$ is {\em monogenic} if $\{1,\theta,\theta^2,\ldots,\theta^{\deg(f)-1}\}$ is a basis for $\Z_K$. 
Thus, from \eqref{Eq:Dis-Dis}, $f(x)$ is monogenic if and only if the index($f$)=1.

  Let $n,a,c\in \Z$, with $ac\ne 0$ and $n\ge 3$. Define 
\begin{align}\label{Def FF}
\begin{split}
   \FF_{n,a,c}(x):&=x^n+cx^{n-1}+acx^{n-2}+a^2cx^{n-3}+\cdots +a^{n-2}cx+a^{n-1}c\\ 
   &=\frac{x^{n+1}+(c-a)x^n-a^nc}{x-a}.
   \end{split}
   \end{align}
   If $\FF_{n,a,c}(x)$ is irreducible over $\Q$, then it follows from \cite[Theorem 4]{HJLith} that 
    \begin{equation}\label{Del FF}
    \begin{gathered}
    \Delta(\FF_{n,a,c})=(-1)^{(3n^2-3n+2)/2}a^{(n-1)(n-2)}c^{n-1}\frac{\D}{(nc+a)^2},\\ 
     \mbox{where} \quad \D:=(-n)^n(c-a)^{n+1}-a^nc(n+1)^{n+1}.
    \end{gathered}
    \end{equation}
  \begin{rem}
    We point out that the irreducibility of $\FF_{n,1,c}(x)$ was first investigated in \cite{Josh}.
  \end{rem} 
   We define  
    \begin{align}\label{P}
   \begin{split}
    P:&=\{p \ \mbox{prime}: \ p\mid \Delta(\FF_{n,a,c})\},\\
    P_0:&=\{p\in P: \ p\mid ac\},\\
    P_1:&=\{p\in P: \ p\nmid ac, \ p\nmid (n+1)\},\\ 
    P_2:&=\{p\in P: \ p\nmid ac, \ p\mid (n+1)\}.
    \end{split}
   \end{align}   Observe that $P=P_0 \cup P_1 \cup P_2$ and that $P_i\cap P_j=\varnothing$ when $i\ne j$.

Recently, under the assumption that $P_2=\varnothing$, Jakhar, Kalwaniya and Yadav \cite{JKY1} proved the following theorem concerning the prime divisors of the index($\FF_{n,a,c}$). 
\begin{thm}\label{Thm:JKY2}
Suppose that $\FF_{n,a,c}(x)$ is irreducible over $\Q$, and that $P_2=\varnothing$. 
Then a prime $p\in P$ does not divide the {\rm index($\FF_{n,a,c}$)} if and only if $p$ satisfies one of the following conditions:
 \begin{enumerate}[label=(\roman*), font=\normalfont]
 \item \label{JKY2:I1} when $p\mid c$, then $p^2\nmid a^{n-1}c$;
 \item \label{JKY2:I2} when $p\nmid c$, then $p\nmid a$ and $p^2\nmid \D$.
 \end{enumerate}
\end{thm} 
\noindent 
 In this note, we extend Theorem \ref{Thm:JKY2} by providing necessary and sufficient conditions for the monogenicity of $\FF_{n,a,c}(x)$ without the restriction that $P_2=\varnothing$. 
 More precisely, we prove the following theorem. 
 \begin{thm}\label{Thm:Main}
 Recall \eqref{Def FF}, \eqref{Del FF} and \eqref{P}. Suppose that $\FF_{n,a,c}(x)$ is irreducible over $\Q$. Then $\FF_{n,a,c}(x)$ is monogenic if and only if all of the following conditions are satisfied: 
 \begin{enumerate}
 \item \label{I0:a}$\abs{a}=1$;
 \item \label{I0:c} $c$ is squarefree;
   \item \label{I1:Main} every prime $p\in P_1$ is such that $p^2\nmid \D$;
   \item \label{I2:Main} every prime $p\in P_2$ is such that  
   \begin{enumerate}
     \item \label{I2A:Main} if $p=2$, then 
     $(a,c \mmod{4})\in \{(1,2),\ (1,3), \ (-1,1), \ (-1,2)\}$; 
     \item \label{I2B:Main} if $p\ge 3$, then $p^3\nmid (c-a)(nc+a)$.    
   \end{enumerate} 
 \end{enumerate}
  \end{thm}

 \section{Preliminaries}\label{Section:Prelims}
 The first result in this section is a standard tool used to determine the monogenicity of a monic irreducible polynomial. 
\begin{thm}[Dedekind's Index Criterion \cite{Cohen}]\label{Thm:Dedekind}
Let $K=\Q(\theta)$ be a number field, $T(x)\in \Z[x]$ the monic minimal polynomial of $\theta$, and $\Z_K$ the ring of integers of $K$. Let $p$ be a prime number and let $\overline{ * }$ denote reduction of $*$ modulo $p$ (in $\Z$, $\Z[x]$ or $\Z[\theta]$). Let
\[\overline{T}(x)=\prod_{i=1}^s\overline{\gamma_i}(x)^{e_i}\]
be the factorization of $T(x)$ modulo $p$ in $\F_p[x]$, and set
\[g(x)=\prod_{i=1}^s\gamma_i(x),\]
where the $\gamma_i(x)\in \Z[x]$ are arbitrary monic lifts of the $\overline{\gamma_i}(x)$. Let $h(x)\in \Z[x]$ be a monic lift of $\overline{T}(x)/\overline{g}(x)$ and set
\[F(x)=\dfrac{g(x)h(x)-T(x)}{p}\in \Z[x].\]
Then 
\[{\rm index}(T)\not \equiv 0 \pmod{p} \Longleftrightarrow \gcd\Bigl(\overline{g},\overline{h},\overline{F}\Bigr)=1 \mbox{ in } \F_p[x].\]
\end{thm}

The next lemma follows from \cite{Apostol}. 
\begin{lemma}\label{cyclo}
Let $d\ge 1$ be an integer, and let $\Phi_d(x)$ denote the cyclotomic polynomial of index $d$. If $p$ is a prime such that $p\nmid d$, then 
   \[\Phi_d(1)\Phi_d(-1)=\Phi_d(1)\Phi_{2d}(1)\not \equiv 0 \pmod{p}.\]
   \end{lemma}

The following lemma will be useful in the proof of Theorem \ref{Thm:Main}. The proof is straightforward and we omit the details.
\begin{lemma}\label{Lem:FF(a)}
  Let $\FF_{n,a,c}(x)$ be as defined in \eqref{Def FF}. Then 
  \begin{align*}
   \FF_{n,a,c}(a)&=a^{n-1}(nc+a) \quad \mbox{and}\\
   \FF_{n,a,c}(-a)&=\left\{\begin{array}{cl}
   a^n & \mbox{if $2\mid n$}\\
   a^{n-1}(c-a) & \mbox{if $2\nmid n$}.  
   \end{array}\right. 
  \end{align*}
  \end{lemma}

\section{The Proof of Theorem \ref{Thm:Main}} 
\begin{proof}
We use Theorem \ref{Thm:Dedekind} with $T(x):=\FF_{n,a,c}(x)$ and $p$ a prime divisor of $\Delta(\FF_{n,a,c})$ to derive necessary and sufficient conditions for when $p\nmid {\rm index}(\FF_{n,a,c})$. We divide the proof into three parts corresponding respectively to $p\in P_i$, where $i\in \{0,1,2\}$. The first two parts, which address $P_0$ and $P_1$, encompass and refine Theorem \ref{Thm:JKY2}, and correspond to items \eqref{I0:a}, \eqref{I0:c} and \eqref{I1:Main} of this theorem. The third part, which addresses $P_2$, represents the lifting of the restriction that $P_2=\varnothing$ in Theorem \ref{Thm:JKY2}, and corresponds to item \eqref{I2:Main} of this theorem.    
\subsection{$\boldsymbol{p\in P_0}$}\text{}\\ 
If $p\mid a$, then $\overline{T}(x)=x^{n-1}(x+c)$, so that we can let $g(x)=x(x+c)$ and $h(x)=x^{n-2}$ in Theorem \ref{Thm:Dedekind}. Hence,
    \[F(x)=\frac{g(x)h(x)-T(x)}{p}=-\frac{acx^{n-2}+a^2cx^{n-3}+\cdots+a^{n-2}cx+a^{n-1}c}{p},\]
     which implies that $\overline{F}(x)=-\frac{a}{p}cx^{n-2}$. Since $\gcd(\overline{g},\overline{h},\overline{F})=x$, we conclude from Theorem \ref{Thm:Dedekind} that $p\mid \mbox{index($T$)}$ and therefore, $T(x)$ is not monogenic. That is, if $T(x)$ is monogenic, then $\abs{a}=1$, confirming item \eqref{I0:a}. 
     
If $p\mid c$, then $\overline{T}(x)=x^n$, and we can let $g(x)=x$ and $h(x)=x^{n-1}$ in Theorem \ref{Thm:Dedekind} to get
    \[F(x)=\frac{g(x)h(x)-T(x)}{p}=-\frac{cx^{n-1}+acx^{n-2}+a^2cx^{n-3}+\cdots+a^{n-2}cx+a^{n-1}c}{p}.\] Therefore, if $p^2\mid c$, then  $\overline{F}(x)=0$ and $\gcd(\overline{g},\overline{h},\overline{F})=x$. Hence, we conclude from Theorem \ref{Thm:Dedekind} that $p\mid \mbox{index($T$)}$ so that $T(x)$ is not monogenic. It follows that if $T(x)$ is monogenic, then $c$ is squarefree, and item \eqref{I0:c} is established. 
    
    For the remainder of the proof of the theorem, we assume that $\abs{a}=1$ and $c$ is squarefree, since we have shown that these conditions are necessary for the monogenicity of $\FF_{n,a,c}(x)$. Note then, with $p\in P_0$, we have that $p\nmid {\rm index}(\FF_{n,a,c})$ if and only if item \ref{JKY2:I1} of Theorem \ref{Thm:JKY2} is satisfied. 
\subsection{$\boldsymbol{p\in P_1}$}\text{}\\
 Then $p\nmid {\rm index}(\FF_{n,a,c})$ if and only if $p^2\nmid \D$ from item \ref{JKY2:I2} of Theorem \ref{Thm:JKY2}, which establishes item \eqref{I1:Main}.

\subsection{$\boldsymbol{p\in P_2}$}\text{}\\ 
Since $p\nmid ac$ and $p\mid (n+1)$, we see from \eqref{Del FF} that $p\mid \D$ so that $p\mid (c-a)$, which implies that $p\mid (nc+a)$.
Suppose that $p^e\mid \mid (n+1)$, and let $N:=(n+1)/p^e$. We make the following straightforward observations. If $2\mid N$, then $p\ge 3$ and $2\nmid n$. If $2\nmid N$, then either $p=2$ and $2\mid (n+1)$, or $p\ge 3$ and $2\nmid(n+1)$. Thus, 
\begin{align}\label{T mod p combined}
   \nonumber
   T(x)&=\frac{x^{n+1}+(c-a)x^n-a^nc}{x-a}\equiv\frac{x^{n+1}-a^{n+1}}{x-a}\equiv \frac{(x^N-a^N)^{p^e}}{x-a} \pmod{p}\\
   & \equiv \left\{\begin{array}{cl}
   (x-a)^{p^e-1}(x+a)^{p^e}\ds \prod_{\substack{d\mid N\\ d>2}} \Phi_d(x)^{p^e} \pmod{p} & \mbox{if $2\mid N$ and $p\ge 3$}\\[2.5em]
    (x+1)^{2^e-1}\ds \prod_{\substack{d\mid N\\ d>1}} \Phi_d(x)^{2^e} \pmod{2} & \mbox{if $2\nmid N$ and $p=2$}\\[2.5em]
    (x-a)^{p^e-1}\ds \prod_{\substack{d\mid N\\ d>1}} \Phi_{d}(ax)^{p^e} \pmod{p} & \mbox{if $2\nmid N$ and $p\ge 3$.} 
   \end{array}\right.
 \end{align}
We now examine the three cases as indicated in \eqref{T mod p combined}. Since $\abs{a}=1$ and $p\nmid d$, we have from Lemma \ref{cyclo} that 
\begin{equation}\label{Phi}
\Phi_d(a)\Phi_d(-a)\not \equiv 0\pmod{p}
\end{equation} in any of these situations.  

Assume first that $2\mid N$, and suppose that  
   \begin{equation}\label{Tbar Case I}
   \overline{T}(x)=(x-a)^{p^e-1}(x+a)^{p^e}\prod_{i=1}^s\overline{\gamma_i}(x)^{e_ip^e},
    \end{equation} where each $\overline{\gamma_i}(x)$ is irreducible in $\F_p[x]$ and $\overline{\gamma_i}(\pm a)\not \equiv 0 \pmod{p}$ from \eqref{Phi}. 
      Then, in Theorem \ref{Thm:Dedekind}, we can let 
   \begin{equation}\label{g and h Tbar Case I}
   \begin{gathered}
   g(x)=(x-a)(x+a)\prod_{i=1}^s\gamma_i(x) \quad \mbox{and}\\
   h(x)=(x-a)^{p^e-2}(x+a)^{p^e-1}\prod_{i=1}^s\gamma_i(x)^{e_ip^e-1},
   \end{gathered}
   \end{equation}
where the $\gamma_i(x)\in \Z[x]$ are arbitrary monic lifts of the $\overline{\gamma_i}(x)$. Since $p^e\ge 3$, it follows from \eqref{g and h Tbar Case I} that $g(a)=h(a)=0$. Hence, from Lemma \ref{Lem:FF(a)}, we have  
  \begin{equation}\label{Case I pF(a)}
  pF(a)=g(a)h(a)-T(a)=-T(a)=-a^{n-1}(nc+a). 
  \end{equation} Similarly, since $g(-a)=h(-a)=0$ and $2\nmid n$, we have that 
  \begin{equation}\label{Case I pF(-a)}
  pF(-a)=g(-a)h(-a)-T(-a)=-T(-a)=-a^{n-1}(c-a).
  \end{equation} 
   
    Since $2\mid N$, we have that $N>1$, so that $s\ge 1$ in \eqref{Tbar Case I}. Let $\gamma_i(\alpha)\equiv 0\pmod{p}$ for some $i$, and some $\alpha$ in an algebraic closure of $\F_p$. Then there exists some divisor $d>1$ of $N$ such that $\Phi_d(\alpha)\equiv 0 \pmod{p}$. Hence, 
\begin{align}\label{alpha methods}
\begin{split}
\alpha^d\equiv 1 \pmod{p} & \ \Longrightarrow \ \alpha^N\equiv 1 \pmod{p}\\
& \ \Longrightarrow \ \alpha^{n+1}=\left(\alpha^N\right)^{p^e}\equiv 1 \pmod{p^2}\\
& \ \Longrightarrow \ \alpha^{n}\equiv \alpha^{-1} \pmod{p^2}.
\end{split}
\end{align} Recall that $\alpha\ne \pm 1$ since $p\nmid d$. Thus, since $p^e\ge 3$, we conclude from \eqref{g and h Tbar Case I} that $g(\alpha)h(\alpha)\equiv 0 \pmod{p^2}$.  
Therefore, from \eqref{alpha methods}, we have
\begin{align}\label{Case I alpha}
\begin{split}
pF(\alpha)&=g(\alpha)h(\alpha)-T(\alpha)\\
&\equiv -T(\alpha) \pmod{p^2}\\
&\equiv -\frac{\alpha^{n+1}+(c-a)\alpha^n-a^nc}{\alpha-a} \pmod{p^2}\\
&\equiv \frac{a(c-a)}{\alpha} \pmod{p^2}.
\end{split}
\end{align} 

Thus, when $2\mid N$ and $p\ge 3$, we deduce from \eqref{Case I pF(a)}, \eqref{Case I pF(-a)}, \eqref{Case I alpha} and Theorem \ref{Thm:Dedekind} that 
\begin{equation}\label{Case I p odd}
 p\nmid {\rm index}(T)   \ \Longleftrightarrow \ p^2\nmid (c-a) \ \mbox{and} \ p^2\nmid (nc+a) \ \Longleftrightarrow \ p^3\nmid (c-a)(nc+a),
\end{equation}
 since $p\mid (c-a)$ and $p\mid (nc+a)$. 

Assume next that $2\nmid N$ with $p=2$, and suppose that 
\begin{equation}\label{Tbar Case II p=2}
   \overline{T}(x)=(x+1)^{2^e-1}\prod_{i=1}^s\overline{\gamma_i}(x)^{e_i2^e},
    \end{equation} where each $\overline{\gamma_i}(x)$ is irreducible in $\F_p[x]$ and $\overline{\gamma_i}(-1)\not \equiv 0 \pmod{2}$ from \eqref{Phi}. 
      Then, in Theorem \ref{Thm:Dedekind}, we can let 
   \begin{equation}\label{g and h Case II p=2}
   \begin{gathered}
   g(x)=(x+1)\prod_{i=1}^s\gamma_i(x) \quad \mbox{and}\\
   h(x)=(x+1)^{2^e-2}\prod_{i=1}^s\gamma_i(x)^{e_i2^e-1},
   \end{gathered}
   \end{equation}
where the $\gamma_i(x)\in \Z[x]$ are arbitrary monic lifts of the $\overline{\gamma_i}(x)$. 

If $e=1$, then observe from \eqref{g and h Case II p=2} that
 $h(-1)\equiv 1\pmod{2}$ by Lemma \ref{Phi}. If $e\ge 2$, then $g(-1)=h(-1)=0$ from \eqref{g and h Case II p=2}, and $n=2^eN-1\equiv 3\pmod{4}$, so that   
 \begin{align}\label{pF(-1) Case II p=2}
 \begin{split}
  2F(-1)&=g(-1)h(-1)-T(-1)=-T(-1)\\ 
  &=\left\{\begin{array}{cl}
   -(c-1)\equiv -c+1 \pmod{4} & \mbox{if $a=1$}\\
  -(nc-1) \equiv c+1 \pmod{4} & \mbox{if $a=-1$}
\end{array}\right.\\
&\equiv -ac+1 \pmod{4}\\
& \equiv -a(c-a) \pmod{4}.   
  \end{split}
  \end{align} 
  by Lemma \ref{Lem:FF(a)}. 
Furthermore, if $N\ge 3$ and $\overline{\gamma_i}(\alpha)\equiv 0 \pmod{2}$, then we see, following \eqref{alpha methods}, that    
\begin{align}\label{Case II p=2 alpha}
\begin{split}
2F(\alpha)&\equiv \frac{a(c-a)}{\alpha} \pmod{4}. 
\end{split}
\end{align} 

Consequently, when $2\nmid N$ and $p=2$, we conclude from \eqref{pF(-1) Case II p=2}, \eqref{Case II p=2 alpha} and Theorem \ref{Thm:Dedekind} that
\begin{align}\label{p=2}
\begin{split}
2\nmid {\rm index}(T)  \ &\Longleftrightarrow \ 4\nmid (c-a)\\ 
  \ &\Longleftrightarrow \ (a,c \mmod{4})\in \{(1,2),\ (1,3), \ (-1,1), \ (-1,2)\}.
  \end{split}
\end{align} 

Finally, assume that $2\nmid N$ with $p\ge 3$. Since the methods are similar to the previous cases, we provide fewer details in this case. When examining the root $x=a$ of $g(x)$ and $h(x)$, we arrive precisely at the conclusion obtained in \eqref{Case I pF(a)}. If $N\ge 3$, when examining a root $\alpha$ of $\overline{\gamma_i}(x)$, we encounter a minor wrinkle when $a=-1$ since $\Phi_d(-x)=\Phi_{2d}(x)$. In this situation, we must address the two possibilities $\alpha^{n+1}\equiv \pm 1\pmod{p^2}$. Nevertheless, together they yield only the single condition given in \eqref{Case I alpha}. 

Combining \eqref{Case I p odd} and \eqref{p=2} completes the proof of item \eqref{I2:Main}, and hence the proof of the theorem.
\end{proof}



\begin{thebibliography}{99}
\bibitem{Apostol} T. M. Ap\'{o}stol, {\em Resultants of cyclotomic polynomials}, Proc. Amer. Math. Soc. {\bf 24} (1970),
457--462.

\bibitem{Cohen} H. Cohen, \emph{A Course in Computational Algebraic Number Theory}, {Springer-Verlag}, 2000.



\bibitem{Josh} J. Harrington, {\em On the factorization of the trinomials $x^n+cx^{n-1}+d$}, {Int. J. Number Theory} {\bf 8} (2012), no. 6, 1513--1518.

\bibitem{HJLith} J. Harrington and L. Jones, {\em Some New Polynomial Discriminant Formulas}, Lith. Math. J. {\bf 61} (2021), no. 4, 483--490.


\bibitem{JKY1} A. Jakhar, R. Kalwaniya and P. Yadav, {\em On the monogenity and Galois group of certain classes of polynomials}, 
Math. Slovaca {\bf 74}, No. 5, 1147--1154 (2024).



    
   
    










\end{thebibliography}
\end{document}